\documentclass{article}
\usepackage{amssymb}
\usepackage{amsmath}
\usepackage{amsthm}
\usepackage{bbold}

\newcounter{enucount}
\newcounter{enuref}
\newcounter{enualph}
\newcounter{enunumb}
\renewcommand{\theenucount}{(\roman{enucount})}

\renewcommand{\theenunumb}{(\arabic{enunumb})}
\renewcommand{\theenualph}{(\alph{enualph})}

\theoremstyle{plain}
\newtheorem{theorem}{Theorem}
\newtheorem{maintheorem}[theorem]{Main Theorem}
\newtheorem{question}[theorem]{Question}

\newtheorem{cruciallemma}[theorem]{Crucial Lemma}

\newtheorem{mainclaim}[theorem]{Main Claim}

\newtheorem{fact}[theorem]{Fact}

\newtheorem{prop}[theorem]{Proposition}

\newtheorem*{theoremA}{Theorem A}
\newtheorem*{theoremB}{Theorem B}

\theoremstyle{definition}

\theoremstyle{remark}

\renewenvironment{enumerate}{\begin{list}{\rm \theenucount}{\usecounter{enucount}
    \setlength{\labelwidth}{1cm}}}
   {\end{list}}
\newenvironment{alphenumerate}{\begin{list}{\rm \theenualph}{\usecounter{enualph}
    \setlength{\labelwidth}{1cm}}}
   {\end{list}}
\newenvironment{arabenumerate}{\begin{list}{\rm \theenunumb}{\usecounter{enunumb}
    \setlength{\labelwidth}{1cm}}}
   {\end{list}}

\newcommand{\re}{{\upharpoonright}}

\newcommand{\A}{{\cal A}}

\renewcommand{\P}{{\cal P}}

\newcommand{\X}{{\cal X}}

\newcommand{\bb}{{\mathfrak b}}
\newcommand{\cc}{{\mathfrak c}}

\newcommand{\hh}{{\mathfrak h}}

\newcommand{\sss}{{\mathfrak s}}

\newcommand{\AAA}{{\mathfrak A}}

\newcommand{\fin}{{\mathrm{fin}}}

\newcommand{\sub}{\subseteq}
\newcommand{\sem}{\setminus}
\newcommand{\twoom}{2^\omega}
\newcommand{\twolom}{2^{<\omega}}

\newcommand{\omoms}{[\omega]^\omega}
\newcommand{\ha}{\,{}\hat{}\,}
\newcommand{\la}{\langle}
\newcommand{\ra}{\rangle}

\newcommand{\noi}{\noindent}

\title{Base matrices of various heights}

\author{J\"org Brendle\thanks{Partially supported by Grants-in-Aid for Scientific Research
   (C) 18K03398, Japan Society for the Promotion of Science. 
      }  \\ 
   Graduate School of System Informatics \\
   Kobe University \\
   Rokko-dai 1-1, Nada-ku \\
   Kobe 657-8501, Japan \\   
   email: {\sf brendle@kobe-u.ac.jp}}

\begin{document}
\maketitle

\begin{abstract}
\noi A classical theorem of Balcar, Pelant, and Simon says that there is a base matrix of height $\hh$, where $\hh$ is the distributivity
number of $\P (\omega) / \fin$. We show that if the continuum $\cc$ is regular, then there is a base matrix of height $\cc$, and that there are
base matrices of any regular uncountable height $\leq \cc$ in the Cohen and random models. This answers questions of Fischer, Koelbing, and
Wohofsky.
\end{abstract}



\section*{Introduction}

A collection $\AAA = \{ \A_\gamma : \gamma < \vartheta \}$ of mad (maximal almost disjoint) families of subsets of the natural numbers $\omega$
is called a {\em distributivity matrix of height $\vartheta$} if
\begin{itemize}
\item $\A_\delta$ {\em refines} $\A_\gamma$ for $\delta \geq \gamma$, i.e. for all $A \in \A_\delta$ there is $B \in \A_\gamma$ with $A \sub^* B$, and
\item there is no {\em common refinement} of the $\A_\gamma$, i.e. no mad family $\A$ refining all the $\A_\gamma$.
\end{itemize}
$\AAA$ is a {\em base matrix} if it is a distributivity matrix and $\bigcup_{\gamma < \vartheta} \A_\gamma$ is {\em dense } in $\P (\omega ) / \fin$,
i.e. for all $B \in \omoms$ there are $\gamma < \vartheta$ and $A \in \A_\gamma$ with $A \sub^* B$. A distributivity (base) matrix $\AAA$ {\em has cofinal branches} if all maximal branches
of $\AAA$ are cofinal (see Section 1 for a formal definition). The {\em distributivity number} $\hh$ of $\P (\omega) / \fin$ is the least cardinal $\kappa$ such that
$\P (\omega) / \fin$ as a forcing notion is not $\kappa$-distributive; equivalently, it is the least $\kappa$ such that there is a collection $\AAA$  of size $\kappa$ of 
mad families without common refinement. Clearly, a distributivity matrix must have height at least $\hh$, and it is easy to see that there is one of height $\hh$ and none of height 
$> \cc$. Furthermore, if there is a distributivity matrix of height $\vartheta$, then there is one of height $cf (\vartheta)$ so that it suffices to consider regular heights.
A famous theorem of Balcar, Pelant, and Simon~\cite{BPS80} (see also~\cite[Theorem 6.20]{Bl10}) says that there is even a base matrix of height $\hh$. It is natural to ask whether there can
consistently be distributivity (base) matrices of other heights, and in interesting recent work, Fischer, Koelbing, and Wohofsky~\cite{FKWta} proved that it is consistent
that $\hh = \omega_1$ and there is a distributivity matrix with cofinal branches of height $\vartheta \leq \cc$ where $\vartheta > \omega_1$ is regular. We show:

\begin{theoremA}
If $\cc$ is regular then there is a base matrix of height $\cc$.
\end{theoremA}

\begin{theoremB}
In the Cohen and random models, there are base matrices of any regular uncountable height $\leq \cc$.
\end{theoremB}

This answers Questions 7.3 and 7.5 of~\cite{FKWta}. Note that our results are incomparable with the one of the latter work. Their construction does not give a base
matrix (it is not clear whether this can be done) while ours doesn't give cofinal branches -- this is impossible in general. In fact, in the Cohen and random models,
$\hh = \omega_1$ is the only cardinal $\vartheta$ for which there is a distributivity (base) matrix of height $\vartheta$ with cofinal branches (this follows from
Fact~\ref{fact1} below).

\section{Preliminaries}

The {\em Cohen model} ({\em random model}, respectively) is the model obtained by adding at least  $\omega_2$ many Cohen (random, resp.) reals to a model
of the continuum hypothesis CH~\cite{BJ95}.

For $A,B \sub \omega$, we say $A$ is {\em almost contained} in $B$, and write $A \sub^* B$, if $A \sem B$ is finite. $A \subsetneq^* B$ if $A \sub^* B$ and
$B \sem A$ is infinite. For an ordinal $\vartheta_0$, $\{ A_\gamma : \gamma < \vartheta_0 \}$ is a {\em $\sub^*$-decreasing chain} of length $\vartheta_0$ if 
$A_\delta \sub^* A_\gamma$ for all $\gamma < \delta < \vartheta_0$. {\em $\subsetneq^*$-decreasing chains} are defined analogously. For a distributivity
matrix $\AAA = \{ \A_\gamma : \gamma < \vartheta \}$ and an ordinal $\vartheta_0 \leq \vartheta$, $\{ A_\gamma : \gamma < \vartheta_0 \}$ is a {\em branch} in $\AAA$
if it is a $\sub^*$-decreasing chain and $A_\gamma \in \A_\gamma$ for $\gamma < \vartheta_0$. A branch is {\em maximal} if it cannot be properly extended to a longer
branch. A branch is {\em cofinal} if $\vartheta_0 = \vartheta$. Every cofinal branch is maximal, but there may be maximal branches that are not cofinal.

\begin{fact}[Folklore]   \label{fact1} 
There are no $\subsetneq^*$-decreasing chains of length $\omega_2$ in $\P (\omega)$ in the Cohen and random models.
\end{fact}

This is proved by an isomorphism-of-names argument using the homogeneity of the Cohen or random algebra.

For $A,B \in \omoms$, $A$ {\em splits} $B$ if both $A \cap B$ and $B \sem A$ are infinite. $\X \sub \omoms$ is a {\em splitting family} if every $B \in \omoms$
is split by a member of $\X$. The {\em splitting number} $\sss$ is the least size of splitting family. It is well-known that $\hh \leq \sss$ (\cite{Bl10} or~\cite{Ha17}).

\begin{fact}[{Folklore, see~\cite{Bl10}, see also~\cite[Proposition 22.13]{Ha17} for Cohen forcing}]   \label{fact2}
After adding at least $\omega_1$ Cohen or random reals to a model of ZFC, $\sss = \omega_1$. (In fact, the first $\omega_1$ generics are a witness for $\sss$.)
\end{fact}

We will prove:

\begin{maintheorem}   \label{maintheorem}
Assume $\vartheta$ is a regular cardinal and
\begin{itemize}
\item[{\rm (A)}] either there is no $\subsetneq^*$-decreasing chain of length $\vartheta$ in $\P (\omega)$,
\item[{\rm (B)}] or $\sss \leq \vartheta$.
\end{itemize}
Then there is a base matrix of height $\vartheta$.
\end{maintheorem}

Clearly, Theorem A follows from part (B) of the main theorem. (We note, however, that splitting families and $\sss \leq \cc$ are not needed in this case, see the comment
at the beginning of the proof of Main Claim~\ref{mainclaim}.) Theorem B follows from either (A) or (B) in view of Facts~\ref{fact1} and~\ref{fact2}.  Note that part (B)
implies that in many other models of set theory there are base matrices of height $\vartheta$ for any regular $\vartheta$ between $\hh$ and $\cc$, e.g. in the Hechler
model (this satisfies $\sss = \omega_1$ by~\cite{BD85}, see also~\cite{Bl10}), or in {\em any} extension by at least $\omega_1$ Cohen or random reals (Fact~\ref{fact2}).
The former is, and the latter may be (depending on the ground model), a model for the failure of (A). We do not know whether (A) $+ \neg$ (B) is consistent but conjecture
that it is. This clearly implies $\sss \geq \bb^{++}$ where $\bb$ is the unbounding number (which is known to be consistent, see~\cite{BF11}).

\section{Proof of Main Theorem}

By recursion on $\alpha < \cc$, we shall construct sets $\Omega_\gamma \sub \cc$ and families $\A_\gamma = \{ A_{\gamma, \alpha } : \alpha \in \Omega_\gamma \}$,
$\gamma < \vartheta$, such that 
\begin{itemize}
\item[(I)] all $\A_\gamma$ are mad,
\item[(II)] if $\gamma < \delta < \vartheta$ and $\beta \in \Omega_\delta$, then there is $\alpha \leq \beta$ in $\Omega_\gamma$ such that $A_{\delta,\beta} \sub^* A_{\gamma,\alpha}$, and
\item[(III)] for all $B \in \omoms$ there are $\gamma < \vartheta$ and $\alpha \in \Omega_\gamma$ such that $A_{\gamma,\alpha} \sub^* B$.
\end{itemize}
This is clearly sufficient. In case (B), let $\{ S_\zeta : \zeta < \nu \}$ be a splitting family with $\nu \leq \vartheta$. Let $\{ (X_\alpha, \xi_\alpha) :
\alpha < \cc \}$ list all pairs $(X,\xi) \in \omoms \times \vartheta$. At stage $\alpha$ of the construction we will have sets $\{ \Omega_\gamma \cap \alpha : \gamma < \vartheta \}$,
ordinals $\{ \eta_\beta : \beta < \alpha \}$, and families $\{ \{ A_{\gamma,\beta} : \beta \in \Omega_\gamma \cap \alpha \} : \gamma < \vartheta \}$ such that
\begin{itemize}
\item[(i$_\alpha$)] $\A^\alpha_\gamma : = \{ A_{\gamma,\beta} : \beta \in \Omega_\gamma \cap \alpha \}$ is a.d. for $\gamma < \vartheta$,
\item[(ii$_\alpha$)] for all $\beta < \alpha$, the set $\{ \gamma : \beta \in \Omega_\gamma \}$ is the interval of ordinals $[\eta_\beta , \max ( \eta_\beta , \xi_\beta ) ]$ and
\begin{itemize}
   \item for $\gamma \in [ \eta_\beta , \max ( \eta_\beta , \xi_\beta ) ]$, $A_{\gamma,\beta} = A_{\eta_\beta,\beta}$,
   \item for $\gamma < \eta_\beta$, there is $\beta' < \beta$ in $\Omega_\gamma$ such that $A_{\eta_\beta, \beta} \subsetneq^* A_{\gamma, \beta '}$, and
\end{itemize}
\item[(iii$_\alpha$)] $A_{\eta_\beta , \beta} \subsetneq^* X_\beta$ and, in case (B), $A_{\eta_\beta,\beta} \sub^* S_\zeta$ or $A_{\eta_\beta,\beta} \sub^* \omega \sem S_\zeta$ where
   $\zeta$ is minimal such that $S_\zeta$ splits $A_{\gamma, \beta'}$ where $\gamma < \eta_\beta$ and $\beta' \in \Omega_\gamma \cap \beta$ are such that
   $A_{\eta_\beta, \beta} \subsetneq^* A_{\gamma, \beta '}$.
\end{itemize}
Let us first see that this suffices for completing the proof: indeed, (II) and (III)  follow from (ii$_\alpha$) and (iii$_\alpha$), respectively. To see (I), fix $\gamma < \vartheta$
and $Y \in \omoms$. Then there is $\alpha < \cc$ such that $(Y,\gamma) = (X_\alpha, \xi_\alpha)$. So $A_{\max (\eta_\alpha , \xi_\alpha ), \alpha } = A_{\eta_\alpha , \alpha}
\sub^* Y$ by (ii$_{\alpha + 1}$) and (iii$_{\alpha + 1}$) and $A_{\max (\eta_\alpha , \xi_\alpha ), \alpha } \sub^* A_{\gamma, \beta} $ for some $\beta \leq \alpha$ by
(ii$_{\alpha + 1}$). Thus $Y \cap A_{\gamma, \beta}$ is infinite, as required.

Next we notice that for $\alpha = 0$ and for limit $\alpha$, there is nothing to show. Hence it suffices to describe the successor step, that is, the construction at
stage $\alpha + 1$, and to prove that (i$_{\alpha + 1}$) through (iii$_{\alpha + 1}$) still hold. Assume $Y \sub^* X_\alpha \cap A_{\gamma,\beta}$
for some $\gamma < \vartheta$ and $\beta \in \Omega_\gamma \cap \alpha$, and let $\delta$ be such that $\gamma < \delta < \vartheta$. We say that
$Y$ {\em splits at } $\delta$ if
\begin{itemize}
\item for all $\gamma ' $ with $\gamma \leq \gamma ' < \delta$ there is $\beta \in \Omega_{\gamma '} \cap \alpha$ such that $Y \sub^* A_{\gamma ' , \beta }$, and
\item there is no $\beta \in \Omega_\delta  \cap \alpha$ such that $Y \sub^* A_{\delta,\beta}$. 
\end{itemize}
We say $Y$ {\em splits below} $\gamma_0 > \gamma$ if there is $\delta$ with $\gamma < \delta < \gamma_0$ such that $Y$ splits at $\delta$. For infinite $Y \sub X_\alpha$, call
$\A^\alpha_\gamma \re Y$ {\em mad} if $\{ Y \cap A_{\gamma,\beta} : \beta \in \Omega_\gamma \cap \alpha$ and $| Y \cap A_{\gamma, \beta} | = \aleph_0 \}$
is a mad family below $Y$. The following is crucial for our construction.

\begin{cruciallemma}   \label{cruciallemma}
Let $\gamma_0 \leq \vartheta$ be an ordinal and let $Y_0 \sub X_\alpha$ be infinite. Assume
\begin{itemize}
\item[{\bf (mad)}] $\A^\alpha_\gamma \re Y_0$ is mad for all $\gamma < \gamma_0$, and
\item[{\bf (split)}] if $Z \sub^* Y_0 \cap A_{\gamma,\beta} $ for some $\gamma < \gamma_0$ and $\beta \in \Omega_\gamma \cap \alpha$, then $Z$ splits below $\gamma_0$.
\end{itemize}
Then: 
\begin{arabenumerate}
\item $cf(\gamma_0) = \omega$,
\item $\A^\alpha_{\gamma_0} \re Y_0$ is not mad, and, more explicitly,
\item there is an infinite $Z \sub Y_0$ almost disjoint from $\A^\alpha_{\gamma_0}$ such that for all $\gamma < \gamma_0$ there is $\beta \in \Omega_\gamma \cap \alpha$
   such that $Z \sub^* A_{\gamma, \beta}$.
\end{arabenumerate}
\end{cruciallemma}

\begin{proof}
By recursion on $n\in\omega$ we construct infinite sets $( Y^0_s : s \in \twolom )$ and $(Y_s : s \in \twolom )$, as well as ordinals $(\delta^0_s : s \in \twolom )$ and
$(\delta_n : n \in \omega )$ such that
\begin{alphenumerate}
\item $Y_s \sub Y^0_s$ and $Y^0_{s \ha i} \sub Y_s$ for $i \in \{ 0,1 \}$,
\item $\delta_n = \max \{ \delta^0_s : |s| = n \} < \gamma_0$ and $\delta^0_{s \ha i} > \delta_{|s|}$ for $i \in \{ 0, 1 \}$,
\item $Y_s$ splits at $\delta^0_s$ and there are distinct $\beta, \beta ' \in \Omega_{\delta^0_s} \cap \alpha$ such that $Y^0_{s \ha 0} = Y_s \cap A_{\delta^0_s, \beta}$
   and $Y^0_{s \ha 1} = Y_s \cap A_{\delta^0_s, \beta '}$ (in particular, $Y^0_{s \ha 0} \cap Y^0_{s \ha 1}$ is finite), and
\item $Y_{s \ha i} = Y^0_{s \ha i} \cap A_{\delta_{|s|} , \beta}$ for some $\beta \in \Omega_{\delta_{|s|}}$, for $i \in \{ 0,1 \}$.
\end{alphenumerate}
We verify we can carry out the construction. In the basic step $n = 0$ and $s = \la\ra$, by {\bf (mad)}, let $Y_{\la\ra} = Y^0_{\la\ra} : = Y_0 \cap A_{0,\beta}$ for some $\beta
\in \Omega_0 \cap \alpha$ such that this intersection is infinite, and let $\delta_0 = \delta^0_{\la\ra} > 0$ be such that $Y_{\la\ra}$ splits at $\delta_0$. By clause
{\bf (split)} we know that $\delta_0 < \gamma_0$.

Suppose $Y_s^0, Y_s$, and $\delta^0_s$ have been constructed for $|s| = n$ and $\delta_n = \max \{ \delta^0_s: |s| = n \} < \gamma_0$ are such that (a) through (d)
hold. We thus know that $Y_s$ splits at $\delta^0_s$ and, by the definition of splitting and clause {\bf (mad)}, we can find distinct $\beta , \beta' \in \Omega_{\delta^0_s}
\cap \alpha$ such that $Y^0_{s \ha 0}:  = Y_s \cap A_{\delta^0_s, \beta}$ and $Y^0_{s \ha 1} : = Y_s \cap A_{\delta^0_s, \beta '}$ are infinite. Using again {\bf (mad)}, we
see that for $i \in \{ 0,1 \}$ there is $\beta \in \Omega_{\delta_n} \cap \alpha$ such that $Y_{s \ha i} : = Y^0_{s \ha i} \cap A_{\delta_{n} , \beta}$ is infinite.
Next let $\delta^0_{s\ha i} > \delta_n$, $i \in \{ 0,1,\}$, such that $Y_{s \ha i} $ splits at $\delta^0_{s \ha i}$. By {\bf (split)}, $\delta^0_{s \ha i} < \gamma_0$.
Finally, let $\delta_{n+1} : = \max \{ \delta^0_{s \ha i} : |s| = n$ and $i \in \{ 0,1 \} \} < \gamma_0$. This completes the construction.

Let $\delta_\omega = \bigcup_n \delta_n$. Clearly $\delta_\omega \leq \gamma_0$ is a limit ordinal of countable cofinality. Next, for $f \in \twoom$ let $Y_f$ be a 
pseudointersection of the $Y_{f \re n}$, $n \in \omega$. If possible, choose $\beta_f \in \Omega_{\delta_\omega} \cap \alpha $ such that $Y_f \cap A_{\delta_\omega , \beta_f}$ is infinite.
By (c) in this construction and by (ii$_\alpha$), we see that if $f \neq f'$ then $\beta_f \neq \beta_{f'}$. However, $\Omega_{\delta_\omega} \cap \alpha$ has size
strictly less than $\cc$, and therefore there is $f \in \twoom$ for which there is no such $\beta_f$. This implies that $\A^\alpha_{\delta_\omega} \re Y_0$ is not mad and,
by {\bf (mad)}, $\gamma_0 = \delta_\omega$. Furthermore, we may let $Z = Y_f$ for this $f$, and the conclusion of the crucial lemma is established.
\end{proof}

We next show:

\begin{mainclaim}  \label{mainclaim}
There is $\gamma < \vartheta$ such that $\A^\alpha_\gamma \re X_\alpha$ is mad.
\end{mainclaim}

\begin{proof}
Note that in case $\vartheta = \cc$ there is nothing to show because by (ii$_\alpha$) we see that a tail of $\Omega_\gamma \cap \alpha$ is empty,
and therefore so is $\A^\alpha_\gamma$ (in fact, the proof of Theorem A is quite a bit simpler than the general argument: there is no need to list the $\xi_\alpha$,
we may simply let $\xi_\alpha = \alpha$,  $\eta_\alpha $ will always be $\leq \alpha$, and the  splitting family is unnecessary).

Hence assume $\vartheta < \cc$. By way of contradiction, suppose all $\A^\alpha_\gamma \re X_\alpha$ are mad.

We split into two cases. First assume there are $\gamma < \vartheta$, $\beta \in \Omega_\gamma \cap \alpha$ and an infinite $Y \sub X_\alpha \cap A_{\gamma,\beta}$
that does not split below $\vartheta$. This means for all $\delta$ with $\gamma \leq \delta < \vartheta$ there is $\beta \in \Omega_\delta \cap \alpha$ such that
$Y \sub^* A_{\delta,\beta}$. By (ii$_\alpha$), we see that there must be a strictly increasing sequence $(\beta_\varepsilon : \varepsilon < \vartheta )$ of ordinals below
$\alpha$ such that for $\varepsilon ' > \varepsilon$,
\begin{itemize}
\item $\eta_{\beta_{\varepsilon '}} > \max ( \eta_{\beta_\varepsilon} , \xi_{\beta_\varepsilon} )$ and $Y \subsetneq^* A_{\eta_{\beta_{\varepsilon '}} , \beta_{\varepsilon '}}
   \subsetneq A_{\eta_{\beta_{\varepsilon }} , \beta_{\varepsilon }}$.
\end{itemize}
In case (A), this contradicts the initial assumption that there are no $\subsetneq^*$-decreasing chains of length $\vartheta$ in $\P (\omega)$. So assume we are
in case (B). Define a sequence $(\zeta_\varepsilon : \varepsilon < \vartheta )$ of ordinals below $\nu$ such that
\begin{itemize}
\item $\zeta_\varepsilon$ is minimal such that $S_{\zeta_\varepsilon}$ splits all $A_{\eta_{\beta_{\varepsilon '}} , \beta_{\varepsilon '}}$ for $\varepsilon ' < \varepsilon$.
\end{itemize}
Using (iii$_\alpha$), we see that $S_{\zeta_\varepsilon}$ does not split $A_{\eta_{\beta_{\varepsilon }} , \beta_{\varepsilon }}$. Therefore the sequence must be strictly
increasing, which is impossible (and thus contradictory) in case $\nu < \vartheta$. If $\nu = \vartheta$ note that there cannot be any $\zeta$ such that $S_\zeta$ splits
$Y$, contradicting the initial assumption that the $S_\zeta$ form a splitting family.

Therefore there are no such $\gamma$, $\beta$, and $Y$. This implies, however, that the assumptions of the crucial lemma with $\gamma_0 = \vartheta$ and $Y_0 = X_\alpha$ are satisfied while
conclusion (1) clearly fails because $\vartheta$ is an uncountable regular cardinal. This final contradiction establishes the main claim.
\end{proof}

We now let $\eta_\alpha : = \min \{ \gamma : \A^\alpha_\gamma \re X_\alpha$ is not mad$\} < \vartheta$. Choose $Y_0 \sub X_\alpha$ infinite
and almost disjoint from all members of $\A^\alpha_{\eta_\alpha}$. Note that $\A^\alpha_\gamma \re Y_0$ is mad for all $\gamma < \eta_\alpha$.
We construct $A_{\eta_\alpha,\alpha}$ by splitting into cases very much like in the proof of the main claim.

\underline{\sf Case a.} There are $\gamma < \eta_\alpha$, $\beta \in \Omega_\gamma \cap \alpha$, and an infinite $Y \sub Y_0 \cap A_{\gamma,\beta}$ that does not
split below $\eta_\alpha$. Then
\begin{itemize}
\item[$(\star)$] for all $\delta$ with $\gamma \leq \delta < \eta_\alpha$ there is $\beta \in \Omega_\delta \cap \alpha$ such that $Y \sub^* A_{\delta,\beta}$.
\end{itemize}
In this case choose infinite $A_{\eta_\alpha,\alpha} \subsetneq^* Y$. In case (B), we additionally let $A_{\eta_\alpha,\alpha} \sub S_\zeta$ where $\zeta$ is minimal such
that $S_\zeta$ splits $Y$.

\underline{\sf Case b.} Whenever $Y \sub Y_0 \cap A_{\gamma,\beta}$ is infinite for some $\gamma < \eta_\alpha$ and $\beta \in \Omega_\gamma \cap \alpha$,
then $Y$ splits below $\eta_\alpha$. Then the assumptions of the crucial lemma are satisfied with $\gamma_0 = \eta_\alpha$ and $Y_0$, and we see that
$cf (\eta_\alpha) = \omega$ and, by (3), there is $Z \sub Y_0$ such that
\begin{itemize}
\item[$(\star\star)$] for all $\gamma < \eta_\alpha$ there is $\beta \in \Omega_\gamma \cap \alpha$ such that $Z \sub^* A_{\gamma,\beta}$.
\end{itemize}
Then we choose infinite $A_{\eta_\alpha,\alpha} \subsetneq^* Z$. In case (B), we additionally let $A_{\eta_\alpha,\alpha} \sub S_\zeta$ where $\zeta$ is minimal such
that $S_\zeta$ splits $Z$.

Next, for all $\gamma$ with $\eta_\alpha \leq \gamma \leq \max (\eta_\alpha,\xi_\alpha)$, we let $A_{\gamma,\alpha} = A_{\eta_\alpha, \alpha}$. Also put
\[ \Omega_\gamma \cap ( \alpha + 1) = \left\{  \begin{array}{ll} \Omega_\gamma \cap \alpha & \mbox{ if } \gamma < \eta_\alpha \mbox{ or } \gamma > \max (\eta_\alpha,\xi_\alpha), \\
   (\Omega_\gamma \cap \alpha) \cup \{ \alpha \}  & \mbox{ if } \eta_\alpha \leq \gamma \leq  \max (\eta_\alpha,\xi_\alpha). \end{array} \right. \]
Then clauses (i$_{\alpha+1}$) and (iii$_{\alpha + 1}$) are immediate, and (ii$_{\alpha +1}$) follows from $(\star)$ or $(\star\star)$, depending on whether we are in
Case a or Case b. This completes the proof of the main theorem. \hfill $\qed$

\section{Further remarks and questions}

Obviously, the main remaining problem is whether the spectrum of heights of base matrices can be non-convex on regular cardinals.

\begin{question}
Is it consistent that for some regular $\vartheta$ with $\hh < \vartheta < \cc$ there is no base (distributivity) matrix of height $\vartheta$?
\end{question}

The simplest instance would be $\hh = \omega_1$ and $\cc = \omega_3$ with no base (distributivity) matrix of height $\omega_2$. By (B) in
Main Theorem~\ref{maintheorem}, this would imply $\sss = \omega_3$.

The proof of Main Theorem~\ref{maintheorem} may look a little like cheating because we do not refine our mad families everywhere when going to
the next level. Thus let us say $\AAA = \{ \A_\gamma : \gamma < \vartheta \}$ is a {\em strict} base (distributivity) matrix if it is a base (distributivity) matrix and for
any $\gamma < \delta < \vartheta$ and any $A \in \A_\delta$ there is $B \in \A_\gamma$ with $A \subsetneq^*  B$. We leave it to the reader to verify the
details of the corresponding versions of Theorems A and B:

\begin{prop}   \label{prop1}
If $\cc \leq \omega_2$ then there is a strict base matrix of height $\cc$.
\end{prop}

\begin{prop}   \label{prop2}
Let $\vartheta$ be a regular uncountable cardinal. In the Cohen and random models, the following are equivalent:
\begin{enumerate}
\item $\vartheta \in \{ \omega_1 , \omega_2 \}$,
\item there is a strict base matrix of height $\vartheta$,
\item there is a strict distributivity matrix of height $\vartheta$.
\end{enumerate}
\end{prop}

To see, e.g., Proposition~\ref{prop2}, note that by Fact~\ref{fact1}, there cannot be a strict distributivity matrix of height $\geq \omega_3$ in either model.
On the other hand, the original construction of~\cite{BPS80} gives a strict base matrix of height $\hh = \omega_1$, and for $\vartheta = \omega_2$ modify the proof
of Main Theorem~\ref{maintheorem} by attaching a $\subsetneq^*$-decreasing chain of length $\zeta_\beta + 1$ to the set $\{ \gamma : \beta \in \Omega_\gamma \} =
[ \eta_\beta , \max (\eta_\beta, \xi_\beta) ]$ where $\eta_\beta + \zeta_\beta = \max (\eta_\beta , \xi_\beta)$. This is possible because such chains exist of any length
$< \omega_2$ in ZFC and $\zeta_\beta < \omega_2$.



\end{document}